\documentclass[12pt]{article}
\usepackage{amsfonts}
\usepackage{latexsym}

\newcommand{\be}{\begin{equation}} \newcommand{\ee}{\end{equation}}
\newcommand{\bea}{\begin{eqnarray}} \newcommand{\eea}{\end{eqnarray}}
\newcommand{\bean}{\begin{eqnarray*}}
  \newcommand{\eean}{\end{eqnarray*}}
\newcommand{\brray}{\begin{array}} \newcommand{\erray}{\end{array}}
\newcommand{\ben}{\begin{equation}{nonumber}}
  \newcommand{\een}{\end{equation}{nonumber}}
\newcommand{\newsection}[1]{\setcounter{equation}{0}
  \setcounter{dfn}{0}
\section{#1}}

\newtheorem{dfn}{Definition}[section] \newtheorem{thm}[dfn]{Theorem}
\newtheorem{lmma}[dfn]{Lemma} \newtheorem{ppsn}[dfn]{Proposition}
\newtheorem{crlre}[dfn]{Corollary} \newtheorem{xmpl}[dfn]{Example}
\newtheorem{rmrk}[dfn]{Remark}

\newcommand{\bdfn}{\begin{dfn}} \newcommand{\bthm}{\begin{thm}}
    \newcommand{\blmma}{\begin{lmma}}
      \newcommand{\bppsn}{\begin{ppsn}}
        \newcommand{\bcrlre}{\begin{crlre}}
          \newcommand{\bxmpl}{\begin{xmpl}}
            \newcommand{\brmrk}{\begin{rmrk}}

              \newcommand{\edfn}{\end{dfn}}
            \newcommand{\ethm}{\end{thm}}
          \newcommand{\elmma}{\end{lmma}}
        \newcommand{\eppsn}{\end{ppsn}}
      \newcommand{\ecrlre}{\end{crlre}}
    \newcommand{\exmpl}{\end{xmpl}} \newcommand{\ermrk}{\end{rmrk}}

\newcommand{\IC}{{\Bbb C}}

 \newcommand{\IN}{{\Bbb N}}
 
 \newcommand{\IR}{{\Bbb R}}
 \newcommand{\IT}{{\Bbb T}}

 \newcommand{\IZ}{{\Bbb Z}}

 \newcommand{\al}{\alpha}

 \newcommand{\del}{\partial}

\newcommand{\cla}{{\cal A}} \newcommand{\clb}{{\cal B}}

\newcommand{\AAZ}{{\cal A}_{\hbar{}}^{\infty}}

  \def 
\bbE
{\mbox{\boldmath $E$}}                     

  \def 
\bbe
{\mbox{\boldmath $e$}}

 \def\a*{{\cal A}_{h,*}} \def\B{{\cal B}(h)}
\def\B1{{\cal B}_1(h)} \def\b{{\cal B}^{s. a. }(h)} \def\b1{{\cal
    B}^{s. a. }_1(h)}

\newcommand{\raro}{\rightarrow}

  \def \qed 
{
  \mbox{}\hfill $\Box$\vspace{1ex}} 

\begin{document}

\author{{\large \sc Partha Sarathi Chakraborty}\\
         Indian Statistical
Institute,\\[-.5ex]
         203, B. T. Road,
Calcutta--700\,035, INDIA\\[-.5ex]
         email:
parthasc$\underline{\mbox{\hspace{.5em}}}$r@isical.ac.in\\[1ex]}

\title{Metrics on  the  Quantum Heisenberg Manifold}
\maketitle

\section {Introduction}


Let G be the Heisenberg group. \\
$G =  \{  \left(  \matrix { 1 & x & z \cr 0 & 1 & y \cr 0 & 0 & 1 \cr } \right)  | x,y,z \in \IR \} $ \\
For a positive integer c, let $H_c$ be the subgroup of $G$ obtained when $ x,y,cz$ are integers. The Heisenberg manifold $M_c$ is the quotient $ G/H_c$. Nonzero Poisson brackets on $M_c$ invariant under left translation by $G$ are parametrized by two real parameters $ \mu,\nu$ with $ \mu^2 + \nu^2 \ne 0$ \cite {RI1}. For each positive integer $c$ and real numbers $\mu,\nu$ Rieffel constructed a C*-algebra $A_{\mu,\nu}^{c,\hbar}$, quantum Heisenberg manifold (QHM) as example of deformation quantization along a Poisson bracket \cite {RI1}. These algebras have further been studied by \cite {AB1}  \cite {AB2} \cite {CS} \cite {W1}. Recently in a series of papers   \cite {RI2}     \cite {RI3}     \cite {RI4}    Rieffel has introduced the notion of compact quantum metric space (CQMS). He has also constructed examples in the case of C*-dynamical systems where the dynamics is driven by a compact Lie group in ergodic way. Now it is also known that Heisenberg group acts ergodically on QHM.
Using this action Weaver attempted to produce examples of CQMS out of QHM. His construction does not completely achieve the goal. Here essentially using the technique of Rieffel in a modified way we construct examples of CQMS from QHM.

Organization of the paper is as follows. In the next section we briefly recall the notion of QHM and the group action. Then  on a suitable dense *-subalgebra we give a *-algebra norm stronger than the C*-norm. In section 3 we recall the definition of CQMS and construct examples out of QHM using the group action and the previously introduced norm.

\newsection{The Quantum Heisenberg Algebra}
Notation: for $x \in \IR$, $e(x)$ stands for $e^{2 \pi  ix}$
\bdfn
For any positive integer c let $S^c$ denote the space of  $ C^{\infty}$  functions
$\Phi : \IR \times \IT \times \IZ  \raro C$ such that \\
a) $\Phi (x+k,y,p)=e(ckpy) \Phi(x,y,p) $ for all $k \in \IZ$ \\
b) for every polynomial $P$ on $\IZ$ and every partial differential operator \\$ \widetilde{X}=\frac{\del^{m+n}}{\del x^m \del y^n} $ on $\IR \times \IT$ the function $ P(p)(\widetilde{X} \Phi) (x,y,p)$ is bounded on $ K \times \IZ$ for any compact subset $K$  of $\IR \times \IT$.\\
For each  $\hbar,\mu,\nu  \in \IR,\mu^2 + \nu^2 \ne 0$, let ${\cla}^{\infty}_{\hbar} $ denote $S^c$ with product and involution defined by
\bea
\label{1}
(\Phi \star \Psi)(x,y,p)= \sum_q \Phi(x-\hbar (q-p) \mu ,y-\hbar (q-p) \nu,q) \Psi (x- \hbar q \mu,y-\hbar q \nu,p-q) \eea
\bea \label{2} \Phi^*(x,y,p)= \bar{\Phi}(x,y,-p) \eea
$\pi : {\cla}^{\infty}_{\hbar}\raro \clb(L^2( \IR \times \IT \times \IZ)) $ given by
\bea \label{3}
(\pi (\Phi) \xi)(x,y,p)= \sum_q \Phi ( x - \hbar ( q-2p)\mu,y- \hbar (q-2p) \nu,q) \xi ( x,y,p-q)
\eea gives a faithful representation of the involutive  algebra $ \AAZ$. \\
${\cla}^{c,\hbar}_{\mu,\nu}=$ norm closure of $\pi( \AAZ )$ is called the Quantum Heisenberg Manifold.\\
$N_{\hbar}$= weak closure of $\pi ( \AAZ )$
\edfn
We will identify $\AAZ$ with $\pi( \AAZ) $ without any mention.\\
Since we are  going to work with fixed parameters $ c, \mu,\nu, \hbar $ we will drop them altogether and denote ${\cla}^{c,\hbar}_{\mu,\nu}$ simply by $\cla_\hbar$ here the subscript remains merely as a reminiscent of  Heisenberg only to distinguish it from a general algebra.\\


 {\bf Action of the heisenberg group:}
For $ \Phi \in S^c, (r,s,t) \in \IR^3 \equiv G$,  (as a topological space)
\bea
\label{4}
(L_{(r,s,t)} \phi )(x,y,p)=e(p(t+cs(x-r)))\phi(x-r,y-s,p) \eea
extends to an ergodic action of the Heisenberg group on $ {\cla}^{c,\hbar}_{\mu,\nu}$.\\


{\bf The Trace:}
 $\tau : \AAZ \raro \IC$, given by $\tau (\phi)= \int^1_0 \int_{\IT} \phi (x,y,0) dx dy $ extends to a faithful normal tracial state on $N_\hbar$.\\
$\tau$ is invariant under the Heisenberg group action.


\bdfn
Let $ \phi \in S^c$, then $ {\| \cdot \|}_{\infty,\infty,1} $ is the norm defined by
$$ {\| \phi  \|}_{\infty,\infty,1}= \sum_{p \in \IZ} \sup_{ x \in \IR, y \in \IT } | \phi (x,y,p) | $$
\edfn

\bppsn
$ {\| \cdot \|}_{\infty,\infty,1} $ is a *-algebra norm on $S^c$.
\eppsn
{\it Proof:} For $ \phi \in S^c$, clearly $ {\| \phi  \|}_{\infty,\infty,1}= {\| \phi^*  \|}_{\infty,\infty,1}=$.
Let $ \phi,\psi  \in S^c$ and $ \phi^\prime  (p) = \sup_{x \in \IR, y \in \IT} | \phi (x,y,p)|,
\psi^\prime  (p) = \sup_{x \in \IR, y \in \IT} | \psi (x,y,p)|$ for $ p \in \IZ$
\bean
& &   |(\Phi \star \Psi)(x,y,p) | \cr
& \le & \sum_q |\Phi(x-\hbar (q-p) \mu ,y-\hbar (q-p) \nu,q)| \times | \Psi (x- \hbar q \mu,y-\hbar q \nu,p-q) | \cr
 & \le & \sum_q  \phi^\prime  (q) \psi^\prime  (p-q) \cr
\eean
Therefore,
\bean
 \sum_p sup_{x \in \IR, y \in \IT}  |(\Phi \star \Psi)(x,y,p) | & \le  &  \sum_p \sum_q \phi^\prime  (q) \psi^\prime  (p-q)  \cr    & =  &  {\| \phi  \|}_{\infty,\infty,1}.{\| \psi  \|}_{\infty,\infty,1} \cr
\eean
 This proves that $ {\| \cdot \|}_{\infty,\infty,1} $ is an algebra norm. \qed.

\bppsn
The topology given by  $ {\| \cdot \|}_{\infty,\infty,1} $ is stronger than the topology given by the C*-norm coming from $\cla_\hbar$
\eppsn
{\it Proof:} It suffices to show for $ \phi \in S^c, \| \phi \| \le {\| \phi  \|}_{\infty,\infty,1}$

Let $\phi^{\prime} : \IZ \raro \IR_+ $ be $ \phi^{\prime} (n)= \sup_{x \in \IR,y \in \IT}|\phi(x,y,n)|$.\\ Then
for $ \xi \in L^2 ( \IR \times \IT \times \IZ )$ \\
$| (\phi \xi)(x,y,p)| \le ( \phi^{\prime} \star | \xi (x,y, . )|)(p)$,\\
where $\star$ denotes convolution on $\IZ$ and $|\xi(x,y, .)|$ is the function $ p \mapsto | \xi(x,y,p)|$.
By Young's inequality
\bean
\| (\phi  \xi) ( x,y, .) \|_{l_2} & \le &  \| \phi^{\prime}  \star | \xi(x,y,.)| \|_{l_2} \cr
   &  \le  &   \|\phi^{\prime} \|_{l_1} \| \xi(x,y,.)\|_{l_2}  \cr
\eean
Therefore , $ \|  \phi  \| \le \| \phi \|_{\infty, \infty ,1 }$, since $\| \phi \|_{\infty,\infty,1}= \| \phi^{\prime} \|_{l_1}$ \qed.


\newsection { Compact Quantum Metric Space : The Example on QHM}
We recall some of the definitions from \cite {RI4}
\bdfn
An order unit space is a real partially ordered vector space $A$  with a distinguished element $e$, the order unit  satisfying \\
(i) ( Order Unit property ) For each $ a \in A$ there is an $ r \in \IR$ such that $ a \le r e $. \\
(ii) ( The Archimidean property ) If $ a \in A$ and if $ a \le r e $ for all $ r \in \IR$ with $ r \ge 0$, then $ a \le 0$.\\
\edfn
\brmrk
The following prescription defines a norm on an order unit space.
$$ \| a \| = \inf \{ r \in \IR | -re \le a \le re \} $$
\ermrk
\bdfn
By a state of an order unit space $(A,e)$ we mean a $\mu \in A^\prime$, the dual of $(A, \| \cdot \|)$ such
that $ \mu ( e) =1 = {\| \mu \|}^\prime$. Here $ { \| \cdot \|}^\prime$ stands for the dual norm on  $A^\prime$.
\edfn
\brmrk
States are automatically positive.
\ermrk
\bxmpl
Motivating example of the above concept is the real subspace of selfadjoint elements in a C*-algebra with the  order structure inherited from the C*-algebra.
\exmpl
\bdfn
Let $(A,e)$ be an order unit space. By a Lip norm on $A$ we mean a seminorm $L,$ on $A$ such that \\
(i) For $ a \in A$, we have $L(a)=0$ iff $ a \in \IR e$ \\
(ii) The topology on $S(A)$ coming from the metric \\ $\rho_L (\mu, \nu ) =sup \{ | \mu (a) -  \nu (a) | | L(a) \le 1 \} $ is the $w^*$ topology.
\edfn
\bdfn
A compact quantum metric space is a pair $(A,L)$ consisting of an order unit space $A$ and a Lip norm $L$ defined on it.
\edfn

The following theorem of Rieffel will be of crucial importance.
\bthm
\label {RT}
({\bf Theorem 4.5 of \cite {RI4}}) Let $L$ be a seminorm on the order unit space $A$ such that $L(a)=0 $ iff $a \in \IR e$. Then $\rho_L$ gives $S(A)$ the $w^*$-topology exactly if\\
(i) $(A,L)$ has finite radius, i.e, $ \exists$ some constant $C$ such that \\ $ \rho_L (\mu, \nu ) \le C \mbox   { for all }  \mu, \nu   \in S(A)$.\\
(ii) $\clb_1= \{ a | L(a) \le 1,  \mbox{ and } \| a \| \le 1 \}$ is totally bounded in $A$ for $ \| \cdot \| $.
\ethm

{\bf General Scheme of Construction :--}\\
Let  $ (A,G ,\al)$ be a $ C^*$ dynamical system with $G$ an n dimensional Lie group acting ergodically.
Let $ A^\infty= \{ a \in A | g \mapsto \al_g(a) $ is smooth \}. Then any $X \in Lie(G)$, the Lie algebra of $G$ induces a derivation $ \delta_X : A^\infty \raro  A^\infty$. Let $ X_1, \ldots , X_n$ be a basis of $Lie(G)$. $L(a)= \vee_{i=1}^n {\| \delta_{X_i} (a) \|}_n $, should be a good candidate for a Lip norm. Here $ {\| \cdot \|}_n$ stands for an algebra  norm on $A$ not necessarily the norm coming from the algebra. This is essentially Rieffel's construction the only modification is he considers the case where $ {\| \cdot \|}_n$  is the algebra norm. Here the problem of construction of Lip norms reduces to construction of the norm $ {\| \cdot \|}_n$ such that $L$ so defined becomes a Lip norm.

{\bf Illustration in the context of quantum Heisenberg manifolds:--} \\
Let  $$X_1=\left(\matrix { 0&1&0 \cr 0&0&0 \cr 0& 0&0\cr}\right),
X_2=\left(\matrix { 0&0&0 \cr 0&0&1\cr 0& 0&0\cr}\right),X_3=\left(\matrix { 0&0& 1 \cr 0&0&0 \cr 0& 0&0\cr}\right)$$ be the canonical basis of the Lie algebra of the Heisenberg group. Then the associated derivations are given by
\bean
\delta_1(\phi)(x,y,p)&=& - \frac {\del \phi } {\del x}(x,y,p) \cr
\delta_2 (\phi) (x,y,p)&=& 2 \pi  i c p x \phi(x,y,p) -  \frac { \del \phi} {\del y}(x,y,p)\cr
\delta_3 (\phi)(x,y,p) &=& 2 \pi   i p  \phi(x,y,p) \cr
\eean
 Notation:-- Henceforth $A$ will stand for ${S^c}_{s.a}$
\bppsn
$ L : {S^c}_{s.a}  \raro \IR_+ $ given by $ L(\phi)= \vee_1^3 {\| \delta_i (\phi) \|}_{\infty,\infty,1} $ is a Lip norm.
\eppsn
{\it Proof:} Since the action is ergodic and ${\| \cdot \|}_{\infty,\infty,1} $ is a norm it follows that $ L(\phi)=0$ iff $\phi$ is a constant multiple of identity. By theorem (\ref {RT}), it suffices to show that every sequence $ \phi_n \in \clb_1= \{\phi | L(\phi) \le 1,  \mbox{ and } \| \phi \| \le 1 \}$ admits a subsequence convergent in the norm coming from the C*algebra.  \\
Let
\bean
f_{1,n}(p)  &=& sup_{x \in \IR, y \in \IT} | \frac {\del \phi_n } {\del x}(x,y,p)  | \cr
f_{2,n}(p)  &=& sup_{x \in \IR, y \in \IT} | 2 \pi  i c p x \phi_n (x,y,p) -  \frac { \del \phi_n } {\del y}(x,y,p) | \cr
f_{3,n}(p)  &=& sup_{x \in \IR, y \in \IT} |  2 \pi   i p  \phi_n(x,y,p) | \cr
\eean
$ L( \phi_n) \le 1 $ is equivalent with $ \sum_p f_{i,n}(p) \le 1 $ for $ i=1,2,3$.
$$ sup_{ | x | \le 2, y\in \IT} |  \frac { \del \phi_n } {\del y}(x,y,p) | \le 4 \pi c f_{3,n}(p) + f_{2,n}(p)  \le 1 +
4 \pi c $$
Now by Arzella-Ascoli theorem thereexists $  \phi : \IR \times \IT \times \IZ \raro \IC $ such that for each $ p \in \IZ$,
$ sup_{| x | \le 2, y\in \IT} | \phi_n (x,y,p)  -  \phi(x,y,p) | \raro 0$.
Clearly $\phi$ satisfies the periodicity condition.\\
{\bf Claim:--} $$ \sum_p sup_{x,y} | \phi ( x,y,p)| < \infty $$
{\it Proof of Claim:-- } Suppose not, then for any $N \in \IN, \exists p_1, \ldots, p_k > N$ such that
$ \sum_i sup_{x,y} | \phi ( x,y,p_i)| > 2$. So, one can take $n$ sufficiently large so that
$$ \sum_{|p| \ge N} sup_{x,y} | \phi_n ( x,y,p)|  \ge \sum_i sup_{x,y} | \phi ( x,y,p_i )|  > 3/2$$
On the otherhand note,
\bea
\label {1}
\sum_{|p| \ge N} sup_{x,y} | \phi_n ( x,y,p)|  & = & \sum_{|p| \ge N} \frac {f_{3,n} (p)} {p}  \cr
  & \le & \frac {1}{N} \sum_p f_{3,n} (p) = \frac {1}{N} \cr
\eea
This leads to a contradiction.   \qed.\\
For $ N \in \IN$ let  $ \phi_{ | p | \le N } ( x,y,p) = \cases { \phi (x,y,p) & for $ | p| \le N$ \cr
0 & for $ | p| >  N$ }$\\
Let $ \epsilon >0$ be given. Choose $N$ such that\\
(i) $ {\| \phi - \phi_{ | p | \le N } \|}_{\infty,\infty,1}  \le \epsilon$,  (ii) $ \frac {1} { N } \le \epsilon$.\\
Then by (\ref {1})   $ {\| \phi_n  - \phi_{ n,| p | \le N } \|}_{\infty,\infty,1}  \le \epsilon, \forall n.$
Now choose $m$ such that \\ for $ m \le n,     {\| \phi_{n ,| p | \le N } - \phi_{ | p | \le N } \|}_{\infty,\infty,1}\le \epsilon.$ Therefore for all $n \ge m , $ \\ $ {\| \phi_n  - \phi \|}_{\infty,\infty,1}\le 3 \epsilon $. Now the result follows from Proposition 2.4. \qed\\
\bppsn
For all $ \mu, \nu \in S( A) $
$$ \rho_L ( \mu, \nu )= \sup \{  | \mu (a) -\nu (a)| | L(a) \le 1 \}   \le 6$$
\eppsn
{\it Proof:} Let $ \phi \in S^c$ be such that $ L ( \phi ) \le 1 $.\\
Then $ L ( L_{(0,0,t)} ( \phi )) \le 1 $. Therefore , $ L ( \int_0^1 L_{(0,0,t)} ( \phi ) dt )  \le 1 $.\\
i.e, $L ( {\phi}^{(0)} ) \le 1$ where ${\phi}^{(0)}(x,y,p) = \delta_{p0} \phi(x,y,p)$.
Recall, $$ \int_0^1 \int_{\IT}  L_{(r,s,0)} ( {\phi}^{(0)}) dr ds ) = \tau ({\phi}^{(0)}) I $$
Let $ f_{3}(p)=   |2 \pi  p \phi (x,y,p) |$, then $ \sum_p  f_{3}(p) \le 1$ since $ L(\phi ) \le 1 $. Now,\\
(i) $ \| \phi - {\phi}^{(0)}\| \le {\| \phi - {\phi}^{(0)}\|}_{\infty,\infty,1} \le \sum_{p \ne 0} \frac {f_3(p)}
{ 2 \pi |p|} \le \sum_p f_3(p) \le 1 $\\
(ii) ${\|  {\phi}^{(0)}-L_{(r,s,0)} {\phi}^{(0)} \|}_{\infty,\infty,1} \le 2$\\
Using these two we get,
\bean
| \mu ( \phi) - \tau ( {\phi}^{(0)}) | & \le &   | \mu ( \phi) - \mu ( {\phi}^{(0)}) |  + | \mu ( \phi) - \mu ( \tau({\phi}^{(0)}) I ) |  \cr
  & \le & \| \phi - {\phi}^{(0)}\|  + \int_0^1 \int_0^1 | \mu ( {\phi}^{(0)}) - \mu ( L_{(r,s,0)} ( {\phi}^{(0)})
| dr ds \cr
& \le & 3. \cr
\eean
This completes the proof. \qed.
\bthm
$((A,I),L)$ is a compact quantum metric space
\ethm
{\it Proof:} Follows from the previous two propositions. \qed.

\end {document}